\documentstyle[amssymb,amsmath,11pt]{article}

\input{tcilatex}

\begin{document}

\title{Generalized Confidence Interval for the Common Coefficient of
Variation }
\author{J. Behboodian* and A. A. Jafari** \\
*Department of Mathematics, Shiraz Islamic Azad University, Shiraz, IRAN\\
{\it email:} Behboodian@stat.susc.ac.ir\\
*Department of Statistics, Shiraz University, Shiraz 71454, IRAN}
\date{}
\maketitle

\begin{abstract}
In this article, we consider the problem of constructing the confidence
interval and testing hypothesis for the common coefficient of variation (CV)
of several normal populations. A new method is suggested using the concepts
of generalized {\it p}-value and generalized confidence interval. Using this
new method and a method proposed by Tian (2005), we obtain a shorter
confidence interval for the common CV. This combination method has good
properties in terms of length and coverage probability compared to other
methods. A simulation study is performed to illustrate properties. Finally,
these methods are applied to two real data sets in medicine.
\end{abstract}


\section{Introduction}

The coefficient of variation (CV) of a random variable $X$, with mean $\mu
\neq 0$ and standard deviation $\sigma $, is defined by the ratio $\dfrac{%
\sigma }{\mu }$. This ratio is an important measure of variation and it is
useful in medicine, biology, physics, finance, toxicology, business,
engineering, and survival analysis, because it is free from the unit of
measurement and it can be used for comparing the variability of two
different populations.

There are different methods for making inferences about the coefficient of
variation. Lehmann (1996) proposed an exact method for a confidence interval
of CV. Vangel (1996), and Wong and Wu (2002) obtained approximate confidence
intervals. Verrill (2003) reviewed the exact approach that is appropriate
for normally distributed data.

Let $X_{ij},$ $i=1,...,k,$ $j=1,...,n_{i},$ be independent normal random
variables with means $\mu _{i}$ and variance $\sigma _{i}^{2}$. Denote the
CV of the $i$th population by $\varphi _{i}=\dfrac{\sigma _{i}}{\mu _{i}}.$
Consider the hypothesis $H_{\circ }:\varphi _{1}=\varphi _{2}=...=\varphi
_{k}$. Miller (1991) proposed an asymptotic test statistic for the $H_{\circ
}$. Fung and Tsang (1998) reviewed several parametric and nonparametric
tests for the equality of CV in $k$ populations. Pardo and Pardo (2000)
introduced a class of test statistics based on R\'{e}nyi's divergence for
this problem. Nairy and Rao (2003) proposed three new tests based on the
inverse sample CV, i.e. $\bar{X}/S$, and discussed about the size and power
comparison of eight tests.

The assumption of the equality of CV's is common in biological and
agricultural experiments (See Fung and Tsang, 1998). Feltz and Miller (1996)
presented one reasonable estimate for the common CV. Ahmed (2002) proposed
six asymptotic estimators for the common CV and discussed on the risk
behavior of the estimators. The generalized {\it p}-value concept was
introduced by Tsui and Weerahandi (1989) and the generalized confidence
interval by Weerahandi (1993). By using these concepts, Tian (2005) proposed
a generalized {\it p}-value and a generalized confidence interval for the
common CV. Verrill and Johnson (2007) obtained confidence bounds on the
common CV and a ratio of two CV's.

In this article, we propose new methods for making inferences about the
common CV. In Section 2, we look at the concepts of generalized {\it p}%
-value and generalized confidence interval. In Section 3, we will first
review the method of Tian (2005) and Verrill and Johnson (2007), briefly,
and then a new method is given to construct a confidence interval and
hypothesis testing for the common CV by using the concept of generalized
variable. Then by combining this new method and the proposed method by Tian
(2005), we obtain a confidence interval for the common CV that has good
properties with respect to other methods. Section 4 is devoted to a
simulation study, to compare the lengths and coverage probabilities of the
four methods that are given in Section 3. Two real medicine examples are
given in Section 5.

\section{Generalized {\it p}-value and generalized confidence interval}

The concept of generalized {\it p}-value was first introduced by Tsui and
Weerahandi (1989) to deal with some nontrivial statistical testing problems.
These problems involve nuisance parameters in such a fashion that the
derivation of a standard pivot is not possible. See also Weerahandi (1995).

Let ${\bf X}$ be a random variable with density function $f({\bf x}\mid
\zeta )$, where $\zeta =(\theta ,{\bf \eta )}$ is a vector of unknown
parameters, $\theta $ is the parameter of interest, and ${\bf \eta }$ is
possibly a vector of nuisance parameters.

Suppose we have the following hypothesis to test: 
\[
H_{\circ }:\theta \leqslant \theta _{\circ }\text{ \ \ }vs\text{ \ \ }%
H_{1}:\theta >\theta _{\circ }\text{,}
\]%
where $\theta _{\circ }$ is a specified value.

Let ${\bf x}$ be the observed value of random variable ${\bf X}$. $T({\bf X};%
{\bf x},\zeta )$ is said to be a generalized variable if the following three
properties hold:

\noindent(i) For fixed ${\bf x}$ and $\zeta =(\theta _{\circ },{\bf \eta )}$%
, the distribution of $T({\bf X};{\bf x},\zeta )$ is free of the nuisance
parameters ${\bf \eta }$.

\noindent(ii) $t_{obs}=T({\bf x};{\bf x},\zeta )$ does not depend on unknown
parameters.

\noindent(iii) For fixed ${\bf x}$ and ${\bf \eta }$, $P(T({\bf X};{\bf x}%
,\zeta )\geqslant t)$ is either stochastically increasing or decreasing in $%
\theta $ for any given $t.$

If $T({\bf X};{\bf x},\zeta )$ is stochastically increasing in $\theta $,
the generalized {\it p}-value is defined as 
\begin{equation}
p=\sup_{\theta \leqslant \theta _{\circ }}\text{ }P(T({\bf X};{\bf x},\theta
,{\bf \eta })\geq t^{\ast })=P(T({\bf X};{\bf x},\theta _{\circ },{\bf \eta }%
)\geq t^{\ast }),
\end{equation}
where $t^{\ast }=T({\bf x};{\bf x},\theta _{\circ },{\bf \eta }).$

To derive a confidence interval for $\theta$, let $T_{c}({\bf X};{\bf x}%
,\theta ,{\bf \eta })$ satisfies the following conditions:

\noindent(i) The distribution of $T_{c}({\bf X};{\bf x},\theta ,{\bf \eta })$
does not depend on any unknown parameters.

\noindent(ii) The observed value of $T_{c}({\bf X};{\bf x},\theta ,{\bf \eta 
})$ is free of nuisance parameters.

\bigskip

Then, $T_{c}({\bf X};{\bf x},\theta ,{\bf \eta })$ is called a generalized
pivotal variable. Further, if $t_{1}$ and $t_{2}$ are such that 
\begin{equation}
P(t_{1}\leqslant T_{c}({\bf X};{\bf x},\theta ,{\bf \eta })\leqslant
t_{2})=1-\alpha ,
\end{equation}
then, $\Theta=\{\theta :t_{1}\leqslant T_{c}({\bf X};{\bf x},\theta ,{\bf %
\eta })\leqslant t_{2}\}$ gives a $100(1-\alpha )\%$ generalized confidence
interval for $\theta .$ For example, if the value of $T_{c}({\bf X};{\bf x}%
,\theta ,{\bf \eta })$ at ${\bf X}={\bf x}$ is $\theta $, then $\{T_{c}({\bf %
x,}\alpha /2),T_{c}({\bf x,}1-\alpha /2)\}$ is a $(1-\alpha )$ confidence
interval for $\theta $, where $T_{c}({\bf x,}\gamma )$ stands for the $%
\gamma $th quantile of $T_{c}({\bf X};{\bf x},\theta ,{\bf \eta })$.

\section{Inferences for $\protect\varphi $}

Consider $k,$ ($k\geqslant 2$) independent random samples ($%
X_{i1},...,X_{in_{i}}$) from $k$ normal populations with means $\mu _{i}$
and unequal variances $\sigma _{i}^{2}$, $i=1,2,...,k.$ For the $i$th
population, let $\bar{X}_{i}=1/n_{i}\sum\limits_{j=1}^{n_{i}}X_{ij}$ and $%
S_{i}^{2}=1/(n_{i}-1)\sum\limits_{j=1}^{n_{i}}(X_{ij}-\bar{X}_{i})^{2}$ be
the sample mean and sample variance, and let $\bar{x}_{i}$ and $s_{i}^{2}$
be the observed value of the sample mean and sample variance, respectively.

Suppose that 
\begin{equation}
\varphi _{1}=\varphi _{2}=...=\varphi _{k}=\varphi
\end{equation}
where $\varphi _{i}=\dfrac{\sigma _{i}}{\mu _{i}}$ and $\varphi $ is the
common CV parameter.

We are interested in developing a confidence interval and hypothesis test
for the common CV, based on the sufficient statistics $\bar{X}_{i}$ and $%
S_{i}^{2}.$

In this section, we first review the method of Tian(2005) and Verrill and
Johnson (2007) for this problem. A new method is introduced for hypothesis
test and confidence interval, regarding $\varphi ,$ by using the concept of
generalized {\it p}-value and generalized confidence interval. At the end,
by combining this method and the method of Tian (2005), we find a new method
which gives a shorter confidence interval.

\subsection{Method of Tian}

Tian (2005) proposed a generalized pivotal variable of the common CV $%
\varphi ,$ by a weighted average of the generalized pivotal variables of CV
based on individual samples as 
\begin{equation}
T_{1}=T_{1}(\bar{X},S;\bar{x},s,\omega)=\frac{\sum\limits_{i=1}^{k}\dfrac{%
n_{i}-1}{\dfrac{\bar{x}_{i}}{s_{i}}\sqrt{\dfrac{U_{i}}{n_{i}-1}}-\dfrac{Z_{i}%
}{\sqrt{n_{i}}}}}{\sum\limits_{i=1}^{k}(n_{i}-1)}=\frac{\sum\limits_{i=1}^{k}%
\dfrac{n_{i}-1}{\dfrac{\bar{x}_{i}}{s_{i}}\dfrac{S_{i}}{\sigma _{i}}-\dfrac{%
\bar{X}_{i}-\mu _{i}}{\sigma _{i}}}}{\sum\limits_{i=1}^{k}(n_{i}-1)},
\end{equation}
where $U_{i}=\dfrac{(n_{i}-1)S_{i}^{2}}{\sigma _{i}^{2}}\sim \chi _{(n_{i}-1)%
\text{ }}^{2}$ and $Z_{i}=\sqrt{n_{i}}(\dfrac{\bar{X}_{i}-\mu _{i}}{\sigma
_{i}})\sim N(0,1),$ $i=1,2,...,k$, and $\bar{X}=(\bar{X}_{1},...,\bar{X}%
_{k}) $ and $S=(S_{1},...,S_{k})$ with the corresponding observed values $%
\bar{x}$ and $s$, and $\omega=(\varphi,\sigma_1,...,\sigma_k)$.

$T_{1} $ is a generalized pivotal variable for $\varphi $ and can be used to
construct a confidence interval and hypothesis test about $\varphi $.

The $(1-\alpha )$ confidence interval for $\varphi $ is 
\[
\left\{ T_{1}(\bar{x},s,\alpha /2),T_{1}(\bar{x},s,1-\alpha /2)\right\} , 
\]
where $T_{1}(\bar{x},s,\gamma )$ is the $\gamma $th quantile of $T_{1}(\bar{X%
},S;\bar{x},s,\omega).$

Tian (2005) evaluated the coverage properties of this confidence interval by
simulation, and showed that the coverage probabilities are close to nominal
level.

For testing 
\[
H_{\circ }:\varphi \leqslant \varphi _{\circ }\quad vs\ \quad H_{1}:\varphi
>\varphi _{\circ }, 
\]
the generalized {\it p}-value based on (5) is 
\begin{equation}
p=P(T_{1}(\bar{X},S;\bar{x},s,\omega)\leqslant \varphi _{\circ }),
\end{equation}
and for testing the hypothesis 
\[
H_{\circ }:\varphi =\varphi _{\circ }\quad vs\ \quad H_{1}:\varphi \neq
\varphi _{\circ }, 
\]
the generalized {\it p}-value is 
\begin{equation}
p=2\min\left\{ P(T_{1}(\bar{X},S;\bar{x},s,\omega)\leqslant \varphi _{\circ
}),P(T_{1}(\bar{X},S;\bar{x},s,\omega)\geqslant \varphi _{\circ })\right\} .
\end{equation}

\subsection{ Method of Verrill and Johnson}

Under the hypothesis in (4) the log-likelihood function can be written as 
\begin{eqnarray}
\ln L(\theta ) &=&\sum\limits_{i=1}^{k}\left( -n_{i}\ln \sigma
_{i}-\sum\limits_{j=1}^{n_{i}}\left( x_{ij}-\frac{\sigma _{i}}{\varphi }%
\right) ^{2}/(2\sigma _{i}^{2})\right) -\frac{n}{2}\ln 2\pi  \nonumber \\
&=&\sum\limits_{i=1}^{k}\left( -n_{i}\ln \sigma
_{i}-\sum\limits_{j=1}^{n_{i}}\left( (n_{i}-1)S_{i}^{2}+n_{i}(\bar{x}_{i}-%
\frac{\sigma _{i}}{\varphi })^{2}\right) /(2\sigma _{i}^{2})\right) -\frac{n%
}{2}\ln 2\pi ,
\end{eqnarray}
where $\theta ^{T}=\omega =(\varphi ,\sigma _{1},...,\sigma _{k})$.

The Newton estimator of $\theta $ is given by 
\[
\theta _{Newt}=-\left. \left[ \frac{\partial ^{2}\ln L}{\partial \theta
_{l}\partial \theta _{m}}\right] ^{-1}\right| _{\theta _{n,c}}\left. \left( 
\begin{array}{c}
{\partial \ln L/\partial \theta _{1}} \\ 
{\vdots } \\ 
{\partial \ln L/\partial \theta _{k+1}}%
\end{array}
\right) \right| _{\theta _{n,c}}+\theta _{n,c}, 
\]
where $\theta _{n,c}$ is any $\sqrt{n}$ - consistent estimator of $\theta $
(Lehmann, 1996).

Verrill and Johnson (2007) obtained an approximate $(1-\alpha )$ confidence
interval for $\varphi $ as 
\begin{equation}
\hat{\varphi}\pm Z_{\alpha /2}\sqrt{\frac{\hat{\varphi}^{4}+\hat{\varphi}%
^{2}/2}{n}}
\end{equation}
where $\hat{\varphi}$ is the first element of $\theta _{Newt}$ and $%
Z_{\alpha /2}$ is appropriate critical value from a standard normal
distribution.

\subsection{A New Method}

Under the hypothesis in (4), we have $X_{ij}\sim N(\eta \sigma _{i},\sigma
_{i}^{2}),$ $i=1,2,...,k,$ where $\eta =\dfrac{1}{\varphi }.$ We can \ show
that if $\sigma _{i}^{2}$'s are known, then the MLE for $\eta $ is 
\begin{equation}
\hat{\eta}=\frac{\sum\limits_{i=1}^{k}\dfrac{n_{i}}{\sigma _{i}}\bar{X}_{i.}%
}{n},
\end{equation}
where $\hat{\eta}\sim N(\eta ,\dfrac{1}{n})$, and $n=\sum%
\limits_{i=1}^{k}n_{i}$.

\noindent{\bf Remark. }If we use $S_{i}^{2}$ as an estimator for $\sigma
_{i}^{2},$ then a reasonable estimator for $\varphi ,$ is 
\begin{equation}
\hat{\varphi}=\frac{n}{\sum\limits_{i=1}^{k}n_{i}\dfrac{\bar{X}_{i.}}{S_{i}}}%
=\frac{n}{\sum\limits_{i=1}^{k}\dfrac{n_{i}}{\hat{\varphi}_{i}}},
\end{equation}
which is a $\sqrt{n} $ - consistent estimator for $\varphi$.

\bigskip

A generalized pivotal variable for estimating $\sigma _{i}^{2}$ can be
expressed as 
\begin{equation}
R_{i}=\sigma _{i}^{2}\frac{s_{i}^{2}}{S_{i}^{2}}=\frac{(n_{i}-1)s_{i}^{2}}{%
U_{i}},\qquad i=1,2,...,k,
\end{equation}
where $U_{i}=\dfrac{(n_{i}-1)S_{i}^{2}}{\sigma _{i}^{2}}\sim \chi _{(n_{i}-1)%
\text{ }}^{2}$ and $s_{i}^{2}$ is an observed value for $S_{i}^{2}$.

We define a generalized pivotal variable for the common CV, $\varphi$, based
on (10) and (12) as 
\begin{equation}
T_{2}=T_{2}(\bar{X},S;\bar{x},s,\omega)=\frac{n}{\sum\limits_{i=1}^{k}\dfrac{%
n_{i}\bar{x}_{i.}}{\sigma _{i}}\dfrac{S_{i}}{s_{i}}-n(\hat{\eta}-\eta )}=%
\frac{n}{\sum\limits_{i=1}^{k}\dfrac{n_{i}\sqrt{U_{i}}}{\sqrt{n_{i}-1}}%
\dfrac{\bar{x}_{i}}{s_{i}}-\sqrt{n}Z},
\end{equation}
where $\bar{X}=(\bar{X}_{1},...,\bar{X}_{k})$ and $S=(S_{1},...,S_{k})$ with
the corresponding observed values $\bar{x}$ and $s$ and $Z=\sqrt{n}(\hat{\eta%
}-\eta )\sim N(0,1)$.

Since $T_{2}(\bar{X},S;\bar{x},s,\omega)$ satisfies the two conditions (i)
the distribution of $T_{2}(\bar{X},S;\bar{x},s,\omega)$ does not depend on
any unknown parameters (ii) the observed value of $T_{2}(\bar{X},S;\bar{x}%
,s,\omega)$ is free of the nuisance parameters, we can use (13) for
constructing a generalized confidence interval for $\varphi .$

The $(1-\alpha )$ confidence interval for $\varphi $ is 
\[
\left\{ T_{2}(\bar{x},s,\alpha /2),T_{2}(\bar{x},s,1-\alpha /2)\right\} , 
\]
where $T_{2}(\bar{x},s,\gamma )$ is the $\gamma $th quantile of $T_{2}(\bar{X%
},S;\bar{x},s,\omega).$

For testing 
\[
H_{\circ }:\varphi \leqslant \varphi _{\circ }\quad vs\ \quad H_{1}:\varphi
>\varphi _{\circ }, 
\]
we use (13) and define 
\begin{equation}
T_{2}^{^{\prime }}(\bar{X},S;\bar{x},s,\omega)=T_{2}(\bar{X},S;\bar{x}%
,s,\omega)-\varphi .
\end{equation}

The distribution of $T_{2}^{^{\prime }}(\bar{X},S;\bar{x},s,\omega)$ is free
from nuisance parameters, the observed value of $T_{2}^{^{\prime }}(\bar{X}%
,S;\bar{x},s,\omega)$, i.e. $t_{obs}^{^{\prime }}$ is zero, and the
distribution function of $T_{2}^{^{\prime }}(\bar{X},S;\bar{x},s,\omega)$ is
an increasing function with respect to $\varphi .$ Therefore $%
T_{2}^{^{\prime }}(\bar{X},S;\bar{x},s,\omega)$ is a generalized variable
for $\varphi $ and the generalized {\it p}-value is 
\begin{equation}
p=P(T_{2}^{^{\prime }}(\bar{X},S;\bar{x},s,\omega)\leqslant
t_{obs}^{^{\prime }}|\varphi =\varphi _{\circ })=P(T_{2}(\bar{X},S;\bar{x}%
,s,\omega)\leqslant \varphi _{\circ }),
\end{equation}
and for testing the hypothesis 
\[
H_{\circ }:\varphi =\varphi _{\circ }\quad vs\ \quad H_{1}:\varphi \neq
\varphi _{\circ }, 
\]
the generalized {\it p}-value based on (14) is 
\begin{equation}
p=2\min\left\{ P(T(\bar{X},S;\bar{x},s,\omega)\geq \varphi _{\circ }),P(T(%
\bar{X},S;\bar{x},s,\omega)\leq \varphi _{\circ })\right\} .
\end{equation}

\subsection{A Combined Method}

For the generalized pivotal variable of $\varphi ,$ we consider a
combination of the generalized pivotal variables in (5) and (13) as follows: 
\begin{equation}
T_{3}(\bar{X},S;\bar{x},s,\omega)=0.5T_{1}(\bar{X},S;\bar{x}%
,s,\omega)+0.5T_{2}(\bar{X},S;\bar{x},s,\omega).
\end{equation}

Since (i) the distribution of $T_1$ and $T_2$ does not on any unknown
parameters (ii) the observed values of $T_1$ and $T_2$ are equal $\varphi$,
therefore $T_{3}(\bar{X},S;\bar{x},s,\omega)$ is a generalized pivotal
variable for common CV $\varphi ,$ and we can use it to obtain a confidence
interval for $\varphi $ and for testing the hypothesis, we define the
generalized variable as 
\begin{equation}
T_{3}^{^{\prime }}=T_{3}-\varphi .
\end{equation}

\subsection{A Computing Algorithm}

For given $k$ independent sample from normal populations, let $i$th sample
contains $n_{i}$ observations with statistics $\overline{x}_{i}$ and $%
s_{i}^{2}$ .

The generalized confidence intervals for $\varphi $ and the generalized $p$%
-value for testing, based on $T_{h}$'s, $h=1,2,3$ can be computed by the
Monte Carlo simulation (See Weerahandi (1995)). The following steps are
given for the generalized variable $T_{3}$ which they are applicable for the
generalized variables $T_{1}$ and $T_{2}:$

1. generate $U_{i}\sim \chi _{(n_{i}-1)}^{2},$ $i=1,...,k.$

2. generate $Z_{i}\sim N(0\,,\,1),$ $i=1,...,k.$

3. generate $Z\sim N(0\,,\,1).$

4. compute $T_{1}$ and $T_{2}$ in (5) and (13).

5. Calculate $T_{3}=0.5T_{1}+0.5T_{2}$.

6. Repeat steps 1 to 5 for $m$ times and obtain $m$ values of $T_{3}$.

\smallskip

Let $T_{3(p)}$ denote the $100p$th percentile of $T_{3}$'s in step 6. Then $%
[T_{3(\alpha /2)},T_{3(1-\alpha /2)}]$ is a Monte Carlo estimate of $%
1-\alpha $ confidence interval for $\varphi .$

The generalized $p$-value for testing \ $\varphi =\varphi _{0}$ \ $vs$ \ $%
\varphi \neq \varphi _{0}$ is $2\min \left\{ P(T_{3}\geq \varphi
_{0}),P(T_{3}\leq \varphi _{0})\right\} $ and the probability $P(T_{3}\geq
\varphi _{0})$ can be estimated by the proportions of the $T_{3}$'s in step
6 that are greater than or equal to $\varphi _{0}.$ Similarly,\thinspace $%
P(T_{3}\leq \varphi _{0})$ can also be estimated.

\section{Simulation Study}

For comparing the coverage probability of the methods introduced in Section
3;

\noindent I) Method of Tian (2005)

\noindent II) Method of Verrill and Johnson (2007)

\noindent III) A method in (13)

\noindent IV) Combined method in (17)

\noindent a simulation study is performed for $k=3$ populations. The data of
size $n_{i},$ $i=1,2,3,$ were generated from normal distributions with mean $%
\mu _{i}$ and variance $\varphi ^{2}\mu _{i}^{2},$ such that all $k$
populations have common CV $\varphi .$ Using 10000 simulations, coverage (C)
probability and average of length (L) estimated. Also we used the algorithm
in Section 3 by $m=5000$ for obtaining the generalized confidence intervals.
The results are given in Tables 1, 2 and 3.

We observed that

\noindent i) The method in (13) and method in (17) produce comparable
results to method of Tian (2005) and method of Verrilla and Johnson (2007).
Therefore, we must apply the four methods to see which one is the best on
the basis of coverage probability and the length of the interval.

\noindent ii) The coverage probabilities of Tian (2005) and the one obtained
by (17) are close to nominal level.

\noindent iii) In some cases the coverage probabilities of the confidence
intervals constructed by (13) are generally lower than the nominal level
although having a slightly shorter average length for the confidence
intervals.

\noindent iv) The coverage probabilities of the confidence intervals
constructed by the method of Verrill and Johnson (2007) smaller than the
nominal level when the sample sizes are small.

\section{Two Real Examples}

{\bf Example 1. }This is the example used by Tian (2005). Actually Fung and
Tsang (1998) showed that the coefficient of variation for MCV in 1995 is not
significantly different from that of 1996. We are interested in making
inferences about the common coefficient of variation of these data. The
sample size, mean, standard deviation and coefficient of variation for MCV
are 63, 84.13, 3.390, 0.0406 from 1995 survey; and 72, 85.68, 2.946, 0.0346
from 1996 survey. These results are derived and explained in detail in the
above articles.

\newpage

\begin{center}
{\bf Table 4.} The confidence intervals for the common CV

\bigskip 
\begin{tabular}{|c|c|c|}
\hline
method & confidence interval & length \\ \hline
Tian (2005) & (0.0347 , 0.0447) & 0.0100 \\ 
Verrill and Johnson (2007) & (0.0324 , 0.0427) & 0.0103 \\ 
New Method in (13) & (0.0332 , 0.0423) & 0.0091 \\ 
Combined Method in (17) & (0.0333, 0.0425) & 0.0092 \\ \hline
\end{tabular}
\end{center}

The estimate of $\varphi $, by different methods, are : (i) Feltz and Miller
(1996), 0.0374 (ii) new method (11), 0.0372 (iii) MLE, 0.0369.

The 95\% confidence intervals for the common CV based on the four methods
are given in Table 4.

\bigskip

\noindent {\bf Example 2.} The data in Appendix D of Fleming and Harrington
(1991) refer to survival times of patients from four hospitals. These data
and their descriptive statistics are given in Table 5.

\begin{center}
{\bf Table 5.} Data and descriptive statistics for survival times of
patients from four hospitals

\bigskip 
\begin{tabular}{|c|l|c|c|c|}
\hline
& Data & $\bar{x}_{i}$ & $s_{i}^{2}$ & $\hat{\varphi}_{i}$ \\ \hline
Hospital 1 & 176 105 266 227 66 & 168.0 & 6880.5 & 0.4937 \\ 
Hospital 2 & 24 5 155 54 & 59.5 & 4460.3 & 1.1224 \\ 
Hospital 3 & 58 64 15 & 45.7 & 714.3 & 0.5853 \\ 
Hospital 4 & 174 42 305 92 30 82 265 237 208 147 & 154.6 & 8894.7 & 0.6100
\\ \hline
\end{tabular}

\bigskip
\end{center}

Nairy and Rao (2003) tested homogeneity of CV's for the hospitals and they
showed that all tests give the same conclusion of accepting the null
hypothesis. Therefore we have common coefficient of variation for these data.

The estimate of $\varphi $, by different methods, are: (i) Feltz and Miller
(1996), 0.6734 (ii) new method (11), 0.6248 (iii) MLE, 0.6015. The estimate
of $\varphi $ based on (11) is close to MLE.

The 95\% confidence intervals for the common CV based on four methods are
given in Table 6. We observe that the length of the interval based on
combined method is shorter than Tian's method. Also the length of the
interval based on the Verrill and Johnson's method is shorter than other
methods but we showed that the coverage probability of this method for small
sample size, is less than nominal level.

\bigskip

\begin{center}
{\bf Table 6.} The confidence intervals for the common CV

\bigskip 
\begin{tabular}{|c|c|c|}
\hline
method & confidence interval & length \\ \hline
Tian (2005) & (-1.7855 , 3.6561) & 5.4416 \\ 
Verrill and Johnson (2007) & (0.4134 , 1.0613) & 0.6479 \\ 
New Method in (13) & (0.4568 , 1.1759) & 0.7191 \\ 
Combined Method in (17 ) & (-0.5457 , 2.2563) & 2.8020 \\ \hline
\end{tabular}
\end{center}

\bigskip

\noindent{\bf Acknowledgement}: The authors thank the editor and referee for
their helpful comments and suggestions. They are also grateful to Islamic
Azad University, Shiraz Branch, Research Council for the support of this
work.

\newpage

\bigskip

\begin{center}
{\bf Table 1}: Simulated coverage probability ($C$) and average length ($L$)
of $95\%$ two sided confidence interval for $\varphi $ (based on 10000
simulation)

{\small 
\begin{tabular}{c}
\\ \hline
\begin{tabular}{cccccc}
\begin{tabular}{c}
$\varphi =0.05$ \\ 
$\mu _{1},\mu _{2},\mu _{3}$%
\end{tabular}
& 
\begin{tabular}{c}
\quad \\ 
$n_{1},n_{2},n_{3}$%
\end{tabular}
& 
\begin{tabular}{c}
I \\ \hline
\begin{tabular}{cc}
$\ \ C$ \ \ \  & $\ \ L$ \ \ 
\end{tabular}%
\end{tabular}
& 
\begin{tabular}{c}
II \\ \hline
\begin{tabular}{cc}
$\ \ C$ \ \ \  & $\ \ L$ \ \ 
\end{tabular}%
\end{tabular}
& 
\begin{tabular}{c}
III \\ \hline
\begin{tabular}{cc}
$\ \ C$ \ \ \  & $\ \ L$ \ \ 
\end{tabular}%
\end{tabular}
& 
\begin{tabular}{c}
IV \\ \hline
\begin{tabular}{cc}
$\ \ C$ \ \ \  & $\ \ L$ \ \ 
\end{tabular}%
\end{tabular}
\\ \hline
$1,1,1$ & 
\begin{tabular}{c}
$5,5,5$ \\ 
$5,5,10$ \\ 
$5,10,30$ \\ 
$10,10,10$ \\ 
$10,20,20$ \\ 
$10,20,30$ \\ 
$20,20,30$ \\ 
$30,30,30$ \\ 
\quad%
\end{tabular}
& 
\begin{tabular}{ll}
$0.950$ & $0.0679$ \\ 
$0.966$ & $0.0487$ \\ 
$0.946$ & $0.0255$ \\ 
$0.953$ & $0.0331$ \\ 
$0.955$ & $0.0228$ \\ 
$0.952$ & $0.0203$ \\ 
$0.948$ & $0.0185$ \\ 
$0.953$ & $0.0158$ \\ 
& 
\end{tabular}
& 
\begin{tabular}{ll}
$0.920$ & $0.0421$ \\ 
$0.934$ & $0.0336$ \\ 
$0.926$ & $0.0218$ \\ 
$0.939$ & $0.0242$ \\ 
$0.947$ & $0.0218$ \\ 
$0.951$ & $0.0194$ \\ 
$0.950$ & $0.0181$ \\ 
$0.953$ & $0.0152$ \\ 
& 
\end{tabular}
& 
\begin{tabular}{ll}
$0.938$ & $0.0441$ \\ 
$0.931$ & $0.0358$ \\ 
$0.930$ & $0.0223$ \\ 
$0.942$ & $0.0279$ \\ 
$0.957$ & $0.0207$ \\ 
$0.940$ & $0.0188$ \\ 
$0.939$ & $0.0173$ \\ 
$0.959$ & $0.0151$ \\ 
& 
\end{tabular}
& 
\begin{tabular}{ll}
$0.952$ & $0.0529$ \\ 
$0.958$ & $0.0403$ \\ 
$0.948$ & $0.0233$ \\ 
$0.951$ & $0.0295$ \\ 
$0.952$ & $0.0214$ \\ 
$0.953$ & $0.0193$ \\ 
$0.952$ & $0.0177$ \\ 
$0.954$ & $0.0153$ \\ 
& 
\end{tabular}
\\ 
$1,1,2$ & 
\begin{tabular}{c}
$5,5,5$ \\ 
$5,5,10$ \\ 
$5,10,30$ \\ 
$10,10,10$ \\ 
$10,20,20$ \\ 
$10,20,30$ \\ 
$20,20,30$ \\ 
$30,30,30$ \\ 
\quad%
\end{tabular}
& 
\begin{tabular}{ll}
$0.964$ & $0.0687$ \\ 
$0.950$ & $0.0504$ \\ 
$0.953$ & $0.0254$ \\ 
$0.955$ & $\ 0.0333$ \\ 
$0.958$ & $0.0228$ \\ 
$0.947$ & $0.0203$ \\ 
$0.946$ & $0.0185$ \\ 
$0.947$ & $0.0159$ \\ 
& 
\end{tabular}
& 
\begin{tabular}{ll}
$0.933$ & $0.0442$ \\ 
$0.931$ & $0.0380$ \\ 
$0.940$ & $0.0241$ \\ 
$0.941$ & $0.0273$ \\ 
$0.946$ & $0.0201$ \\ 
$0.950$ & $0.0185$ \\ 
$0.954$ & $0.0179$ \\ 
$0.949$ & $0.0150$ \\ 
& 
\end{tabular}
& 
\begin{tabular}{ll}
$0.935$ & $0.0446$ \\ 
$0.928$ & $0.0366$ \\ 
$0.938$ & $0.0223$ \\ 
$0.943$ & $0.0281$ \\ 
$0.948$ & $0.0206$ \\ 
$0.942$ & $0.0188$ \\ 
$0.944$ & $0.0173$ \\ 
$0.947$ & $0.0151$ \\ 
& 
\end{tabular}
& 
\begin{tabular}{ll}
$0.969$ & $0.0535$ \\ 
$0.948$ & $0.0415$ \\ 
$0.949$ & $0.0232$ \\ 
$0.951$ & $0.0297$ \\ 
$0.953$ & $0.0213$ \\ 
$0.948$ & $0.0192$ \\ 
$0.946$ & $0.0176$ \\ 
$0.947$ & $0.0153$ \\ 
& 
\end{tabular}
\\ 
$1,5,10$ & 
\begin{tabular}{c}
$5,5,5$ \\ 
$5,5,10$ \\ 
$5,10,30$ \\ 
$10,10,10$ \\ 
$10,20,20$ \\ 
$10,20,30$ \\ 
$20,20,30$ \\ 
$30,30,30$ \\ 
\quad%
\end{tabular}
& 
\begin{tabular}{ll}
$0.955$ & $0.0685$ \\ 
$0.943$ & $0.0500$ \\ 
$0.955$ & $0.0254$ \\ 
$0.945$ & $0.0325$ \\ 
$0.957$ & $0.0229$ \\ 
$0.960$ & $0.0202$ \\ 
$0.953$ & $0.0185$ \\ 
$0.950$ & $0.0158$ \\ 
& 
\end{tabular}
& 
\begin{tabular}{ll}
$0.910$ & $0.0414$ \\ 
$0.930$ & $0.0401$ \\ 
$0.938$ & $0.0218$ \\ 
$0.941$ & $0.0280$ \\ 
$0.949$ & $0.0218$ \\ 
$0.951$ & $0.0189$ \\ 
$0.953$ & $0.0178$ \\ 
$0.955$ & $0.1570$ \\ 
& 
\end{tabular}
& 
\begin{tabular}{ll}
$0.922$ & $0.0447$ \\ 
$0.922$ & $0.0365$ \\ 
$0.936$ & $0.0222$ \\ 
$0.940$ & $0.0281$ \\ 
$0.941$ & $0.0208$ \\ 
$0.945$ & $0.0187$ \\ 
$0.947$ & $0.0173$ \\ 
$0.939$ & $0.0151$ \\ 
& 
\end{tabular}
& 
\begin{tabular}{ll}
$0.963$ & $0.0535$ \\ 
$0.948$ & $0.0412$ \\ 
$0.952$ & $0.0231$ \\ 
$0.946$ & $0.0298$ \\ 
$0.950$ & $0.0214$ \\ 
$0.958$ & $0.0192$ \\ 
$0.951$ & $0.0177$ \\ 
$0.950$ & $0.0153$ \\ 
& 
\end{tabular}%
\end{tabular}
\\ \hline
\end{tabular}
}
\end{center}

\bigskip

\newpage

\begin{center}
{\bf Table 2}: Simulated coverage probability ($C$) and average length ($L$)
of $95\%$ two sided confidence interval for $\varphi $ (based on 10000
simulation)

{\small 
\begin{tabular}{c}
\\ \hline
\begin{tabular}{cccccc}
\begin{tabular}{c}
$\varphi =0.3$ \\ 
$\mu _{1},\mu _{2},\mu _{3}$%
\end{tabular}
& 
\begin{tabular}{c}
\quad \\ 
$n_{1},n_{2},n_{3}$%
\end{tabular}
& 
\begin{tabular}{c}
I \\ \hline
\begin{tabular}{cc}
$\ \ C$ \ \ \  & $\ \ L$ \ \ 
\end{tabular}%
\end{tabular}
& 
\begin{tabular}{c}
II \\ \hline
\begin{tabular}{cc}
$\ \ C$ \ \ \  & $\ \ L$ \ \ 
\end{tabular}%
\end{tabular}
& 
\begin{tabular}{c}
III \\ \hline
\begin{tabular}{cc}
$\ \ C$ \ \ \  & $\ \ L$ \ \ 
\end{tabular}%
\end{tabular}
& 
\begin{tabular}{c}
IV \\ \hline
\begin{tabular}{cc}
$\ \ C$ \ \ \  & $\ \ L$ \ \ 
\end{tabular}%
\end{tabular}
\\ \hline
$1,1,1$ & 
\begin{tabular}{c}
$5,5,5$ \\ 
$5,5,10$ \\ 
$5,10,30$ \\ 
$10,10,10$ \\ 
$10,20,20$ \\ 
$10,20,30$ \\ 
$20,20,30$ \\ 
$30,30,30$ \\ 
\quad%
\end{tabular}
& 
\begin{tabular}{ll}
$0.967$ & $0.6581$ \\ 
$0.965$ & $0.4333$ \\ 
$0.959$ & $0.1842$ \\ 
$0.955$ & $0.2404$ \\ 
$0.948$ & $0.1581$ \\ 
$0.953$ & $0.1239$ \\ 
$0.955$ & $0.1237$ \\ 
$0.964$ & $0.1057$ \\ 
& 
\end{tabular}
& 
\begin{tabular}{ll}
$0.930$ & $0.3712$ \\ 
$0.932$ & $0.2523$ \\ 
$0.930$ & $0.1371$ \\ 
$0.936$ & $0.1726$ \\ 
$0.948$ & $0.1421$ \\ 
$0.952$ & $0.1226$ \\ 
$0.949$ & $0.1247$ \\ 
$0.953$ & $0.0980$ \\ 
& 
\end{tabular}
& 
\begin{tabular}{ll}
$0.926$ & $0.2956$ \\ 
$0.931$ & $0.2405$ \\ 
$0.932$ & $0.1451$ \\ 
$0.941$ & $0.1831$ \\ 
$0.935$ & $0.1358$ \\ 
$0.946$ & $0.1127$ \\ 
$0.943$ & $0.1125$ \\ 
$0.951$ & $0.0861$ \\ 
& 
\end{tabular}
& 
\begin{tabular}{ll}
$0.956$ & $0.4371$ \\ 
$0.964$ & $0.3091$ \\ 
$0.946$ & $0.1556$ \\ 
$0.954$ & $0.1971$ \\ 
$0.946$ & $0.1383$ \\ 
$0.949$ & $0.1122$ \\ 
$0.948$ & $0.1120$ \\ 
$0.956$ & $0.0972$ \\ 
& 
\end{tabular}
\\ 
$1,1,2$ & 
\begin{tabular}{c}
$5,5,5$ \\ 
$5,5,10$ \\ 
$5,10,30$ \\ 
$10,10,10$ \\ 
$10,20,20$ \\ 
$10,20,30$ \\ 
$20,20,30$ \\ 
$30,30,30$ \\ 
\quad%
\end{tabular}
& 
\begin{tabular}{ll}
$0.968$ & $0.6535$ \\ 
$0.961$ & $0.4285$ \\ 
$0.956$ & $0.1876$ \\ 
$0.960$ & $0.2398$ \\ 
$0.955$ & $0.1582$ \\ 
$0.956$ & $0.1385$ \\ 
$0.947$ & $0.1243$ \\ 
$0.952$ & $0.1059$ \\ 
& 
\end{tabular}
& 
\begin{tabular}{ll}
$0.933$ & $0.3562$ \\ 
$0.938$ & $0.2141$ \\ 
$0.928$ & $0.1252$ \\ 
$0.948$ & $0.1736$ \\ 
$0.953$ & $0.1398$ \\ 
$0.951$ & $0.1147$ \\ 
$0.957$ & $0.1421$ \\ 
$0.953$ & $0.1149$ \\ 
& 
\end{tabular}
& 
\begin{tabular}{ll}
$0.925$ & $0.2942$ \\ 
$0.927$ & $0.2395$ \\ 
$0.933$ & $0.1449$ \\ 
$0.949$ & $0.1831$ \\ 
$0.946$ & $0.1358$ \\ 
$0.944$ & $0.1225$ \\ 
$0.940$ & $0.1130$ \\ 
$0.946$ & $0.0984$ \\ 
& 
\end{tabular}
& 
\begin{tabular}{ll}
$0.966$ & $0.4337$ \\ 
$0.961$ & $0.3072$ \\ 
$0.947$ & $0.1558$ \\ 
$0.958$ & $0.1969$ \\ 
$0.947$ & $0.1386$ \\ 
$0.947$ & $0.1232$ \\ 
$0.945$ & $0.1225$ \\ 
$0.948$ & $0.1037$ \\ 
& 
\end{tabular}
\\ 
$1,5,10$ & 
\begin{tabular}{c}
$5,5,5$ \\ 
$5,5,10$ \\ 
$5,10,30$ \\ 
$10,10,10$ \\ 
$10,20,20$ \\ 
$10,20,30$ \\ 
$20,20,30$ \\ 
$30,30,30$ \\ 
\quad%
\end{tabular}
& 
\begin{tabular}{ll}
$0.961$ & $0.6462$ \\ 
$0.966$ & $0.4319$ \\ 
$0.948$ & $0.1910$ \\ 
$0.956$ & $0.2433$ \\ 
$0.949$ & $0.1569$ \\ 
$0.952$ & $0.1381$ \\ 
$0.952$ & $0.1237$ \\ 
$0.955$ & $0.1064$ \\ 
& 
\end{tabular}
& 
\begin{tabular}{ll}
$0.914$ & $0.2301$ \\ 
$0.922$ & $0.2415$ \\ 
$0.928$ & $0.1453$ \\ 
$0.936$ & $0.1722$ \\ 
$0.943$ & $0.1251$ \\ 
$0.952$ & $0.1326$ \\ 
$0.960$ & $0.1362$ \\ 
$0.953$ & $0.1106$ \\ 
& 
\end{tabular}
& 
\begin{tabular}{ll}
$0.927$ & $0.2932$ \\ 
$0.928$ & $0.2382$ \\ 
$0.929$ & $0.1462$ \\ 
$0.948$ & $0.1844$ \\ 
$0.945$ & $0.1352$ \\ 
$0.942$ & $0.1220$ \\ 
$0.950$ & $0.1126$ \\ 
$0.948$ & $0.0988$ \\ 
& 
\end{tabular}
& 
\begin{tabular}{ll}
$0.960$ & $0.4305$ \\ 
$0.964$ & $0.3079$ \\ 
$0.948$ & $0.1576$ \\ 
$0.946$ & $0.1986$ \\ 
$0.948$ & $0.1376$ \\ 
$0.952$ & $0.1357$ \\ 
$0.953$ & $0.1230$ \\ 
$0.956$ & $0.1022$ \\ 
& 
\end{tabular}%
\end{tabular}
\\ \hline
\end{tabular}
}
\end{center}

\bigskip

\newpage

\begin{center}
{\bf Table 3}: Simulated coverage probability ($C$) and average length ($L$)
of $95\%$ two sided confidence interval for $\varphi $ (based on 10000
simulation)

{\small 
\begin{tabular}{c}
\\ \hline
\begin{tabular}{cccccc}
\begin{tabular}{c}
$\varphi =0.5$ \\ 
$\mu _{1},\mu _{2},\mu _{3}$%
\end{tabular}
& 
\begin{tabular}{c}
\quad \\ 
$n_{1},n_{2},n_{3}$%
\end{tabular}
& 
\begin{tabular}{c}
I \\ \hline
\begin{tabular}{cc}
$\ \ C$ \ \ \  & $\ \ L$ \ \ 
\end{tabular}%
\end{tabular}
& 
\begin{tabular}{c}
II \\ \hline
\begin{tabular}{cc}
$\ \ C$ \ \ \  & $\ \ L$ \ \ 
\end{tabular}%
\end{tabular}
& 
\begin{tabular}{c}
III \\ \hline
\begin{tabular}{cc}
$\ \ C$ \ \ \  & $\ \ L$ \ \ 
\end{tabular}%
\end{tabular}
& 
\begin{tabular}{c}
IV \\ \hline
\begin{tabular}{cc}
$\ \ C$ \ \ \  & $\ \ L$ \ \ 
\end{tabular}%
\end{tabular}
\\ \hline
$1,1,1$ & 
\begin{tabular}{c}
$5,5,5$ \\ 
$5,5,10$ \\ 
$5,10,30$ \\ 
$10,10,10$ \\ 
$10,20,20$ \\ 
$10,20,30$ \\ 
$20,20,30$ \\ 
$30,30,30$ \\ 
\quad%
\end{tabular}
& 
\begin{tabular}{ll}
$0.968$ & $2.8818$ \\ 
$0.969$ & $1.5770$ \\ 
$0.967$ & $0.5152$ \\ 
$0.956$ & $0.5907$ \\ 
$0.951$ & $0.2818$ \\ 
$0.954$ & $0.3322$ \\ 
$0.959$ & $0.2517$ \\ 
$0.957$ & $0.2079$ \\ 
& 
\end{tabular}
& 
\begin{tabular}{ll}
$0.921$ & $0.7262$ \\ 
$0.930$ & $0.6221$ \\ 
$0.926$ & $0.2451$ \\ 
$0.942$ & $0.3726$ \\ 
$0.943$ & $0.2471$ \\ 
$0.946$ & $0.2658$ \\ 
$0.956$ & $0.2217$ \\ 
$0.963$ & $0.1923$ \\ 
& 
\end{tabular}
& 
\begin{tabular}{ll}
$0.932$ & $0.5887$ \\ 
$0.936$ & $0.4699$ \\ 
$0.935$ & $0.2775$ \\ 
$0.949$ & $0.2331$ \\ 
$0.948$ & $0.2328$ \\ 
$0.944$ & $0.2585$ \\ 
$0.951$ & $0.2143$ \\ 
$0.947$ & $0.1849$ \\ 
& 
\end{tabular}
& 
\begin{tabular}{ll}
$0.966$ & $1.5547$ \\ 
$0.959$ & $0.9118$ \\ 
$0.954$ & $0.3514$ \\ 
$0.953$ & $0.4213$ \\ 
$0.945$ & $0.2302$ \\ 
$0.952$ & $0.2768$ \\ 
$0.949$ & $0.2088$ \\ 
$0.948$ & $0.1764$ \\ 
& 
\end{tabular}
\\ 
$1,1,2$ & 
\begin{tabular}{c}
$5,5,5$ \\ 
$5,5,10$ \\ 
$5,10,30$ \\ 
$10,10,10$ \\ 
$10,20,20$ \\ 
$10,20,30$ \\ 
$20,20,30$ \\ 
$30,30,30$ \\ 
\quad%
\end{tabular}
& 
\begin{tabular}{ll}
$0.969$ & $2.8874$ \\ 
$0.955$ & $1.5216$ \\ 
$0.967$ & $0.4875$ \\ 
$0.960$ & $0.5884$ \\ 
$0.948$ & $0.3416$ \\ 
$0.949$ & $0.2951$ \\ 
$0.948$ & $0.2487$ \\ 
$0.951$ & $0.2083$ \\ 
& 
\end{tabular}
& 
\begin{tabular}{ll}
$0.925$ & $0.9531$ \\ 
$0.936$ & $0.7216$ \\ 
$0.946$ & $0.2741$ \\ 
$0.952$ & $0.4651$ \\ 
$0.953$ & $0.2821$ \\ 
$0.946$ & $0.2212$ \\ 
$0.956$ & $0.2317$ \\ 
$0.952$ & $0.1851$ \\ 
& 
\end{tabular}
& 
\begin{tabular}{ll}
$0.928$ & $0.5956$ \\ 
$0.925$ & $0.4652$ \\ 
$0.945$ & $0.2744$ \\ 
$0.938$ & $0.3534$ \\ 
$0.950$ & $0.2584$ \\ 
$0.941$ & $0.2328$ \\ 
$0.944$ & $0.2127$ \\ 
$0.948$ & $0.1852$ \\ 
& 
\end{tabular}
& 
\begin{tabular}{ll}
$0.962$ & $1.5651$ \\ 
$0.958$ & $0.8846$ \\ 
$0.958$ & $0.3391$ \\ 
$0.945$ & $0.4170$ \\ 
$0.952$ & $0.2579$ \\ 
$0.953$ & $0.2358$ \\ 
$0.946$ & $0.2017$ \\ 
$0.952$ & $0.1772$ \\ 
& 
\end{tabular}
\\ 
$1,5,10$ & 
\begin{tabular}{c}
$5,5,5$ \\ 
$5,5,10$ \\ 
$5,10,30$ \\ 
$10,10,10$ \\ 
$10,20,20$ \\ 
$10,20,30$ \\ 
$20,20,30$ \\ 
$30,30,30$ \\ 
\quad%
\end{tabular}
& 
\begin{tabular}{ll}
$0.967$ & $2.8634$ \\ 
$0.966$ & $1.6602$ \\ 
$0.962$ & $0.5244$ \\ 
$0.953$ & $0.5979$ \\ 
$0.958$ & $0.3399$ \\ 
$0.953$ & $0.2943$ \\ 
$0.950$ & $0.2515$ \\ 
$0.953$ & $0.2094$ \\ 
& 
\end{tabular}
& 
\begin{tabular}{ll}
$0.930$ & $0.6528$ \\ 
$0.928$ & $0.6211$ \\ 
$0.940$ & $0.3224$ \\ 
$0.937$ & $0.3852$ \\ 
$0.944$ & $0.2634$ \\ 
$0.952$ & $0.2471$ \\ 
$0.957$ & $0.2716$ \\ 
$0.951$ & $0.1928$ \\ 
& 
\end{tabular}
& 
\begin{tabular}{ll}
$0.936$ & $0.5888$ \\ 
$0.926$ & $0.4633$ \\ 
$0.937$ & $0.2772$ \\ 
$0.945$ & $0.3549$ \\ 
$0.942$ & $0.2569$ \\ 
$0.945$ & $0.2337$ \\ 
$0.947$ & $0.2135$ \\ 
$0.952$ & $0.1869$ \\ 
& 
\end{tabular}
& 
\begin{tabular}{ll}
$0.964$ & $1.5513$ \\ 
$0.964$ & $0.9465$ \\ 
$0.959$ & $0.3559$ \\ 
$0.949$ & $0.4231$ \\ 
$0.946$ & $0.2663$ \\ 
$0.951$ & $0.2362$ \\ 
$0.951$ & $0.2084$ \\ 
$0.949$ & $0.1784$ \\ 
& 
\end{tabular}%
\end{tabular}
\\ \hline
\end{tabular}
}
\end{center}

\end{document}